\documentclass[11pt]{article}
\usepackage{amsmath}
\usepackage{latexsym}
\usepackage{amsfonts}
\usepackage{amssymb}

\newcommand{\nf}{\hspace*{-5pt}}

\def \beq { \begin{equation} }
\def \eeq { \end{equation} }

\def\overset#1\to#2{\mathrel{\mathop{#2}\limits^{#1}}}
\def\underset#1\to#2{\mathrel{\mathop{#2}\limits_{#1}}}

\def \varinjlim {\underset{\longrightarrow}\to{lim}}
\def \varprojlim {\underset{\longleftarrow}\to{lim}}

\def \li {\varinjlim}
\def \lp {\varprojlim}

\def \lli#1 {\mathrel{\mathop{\li}\limits_{#1}}}
\def \llp#1 {\mathrel{\mathop{\lp}\limits_{#1}}}

\def \empty {\emptyset}
\def \leqs {\leqslant}
\def \rest {\restriction}

\def \im {\Rightarrow}

\def \lrsa {\leftrightsquigarrow }
\def \rat {\rightarrowtail }

\def \hookr {\hookrightarrow }

\def \rlas {\rightleftarrows}

\newcommand{\hem}{\hspace*{1em}}

\newcommand{\hfl}{\hspace*{\fill}}

\newlength{\dede}     

\newcommand{\hp}{\hspace*{\parindent}} 
\newcommand{\nhp}{\hspace*{-\parindent}}

\newcommand{\N}{\mbox{$ \mathbb N$}}   

\newcommand{\lra}{\mbox{$\longrightarrow$}}  

\newcommand{\sub}{\mbox{$\subseteq$}}

\newcommand{\bcup}{\mbox{$ \bigcup$}}

\def \empty {\emptyset}

\newcommand{\cI}{\mbox{$\cal I$}}

\newcommand{\cL}{\mbox{$\cal L$}}

\newcommand{\cA}{\mbox{$\cal A$}}

\newcommand{\cQ}{\mbox{$\cal Q$}}

\newcommand{\cS}{\mbox{$\cal S$}}
\newcommand{\cT}{\mbox{$\cal T$}}

\newcommand{\all}{\mbox{$\forall$}}

\newcommand{\Lra}{\mbox{$\Leftrightarrow$}}

\def \rest {\restriction}
\def \ra {\rightarrow}

\def \im {\Rightarrow}

\def \lrsa {\leftrightsquigarrow }
\def \rat {\rightarrowtail }

\def \hookr {\hookrightarrow }

\def\overset#1\to#2{\mathrel{\mathop{#2}\limits^{#1}}}
\def\underset#1\to#2{\mathrel{\mathop{#2}\limits_{#1}}}

\newcommand{\und}[1]{\raisebox{-.2ex}{\underline{\raisebox{.2ex}{#1}}}}  

\def \leqs {\leqslant}
\def \geqs {\geqslant}

\def \varinjlim {\underset{\longrightarrow}\to{lim}}
\def \varprojlim {\underset{\longleftarrow}\to{lim}}

\def \li {\varinjlim}
\def \lp {\varprojlim}

\def \lli#1 {\mathrel{\mathop{\li}\limits_{#1}}}
\def \llp#1 {\mathrel{\mathop{\lp}\limits_{#1}}}

\newcommand{\td}[1]{\mbox{$\widetilde{#1}$}}

\newtheorem{Th}{Theorem}[section] 
\newtheorem{Co}[Th]{Corollary}
\newtheorem{Df}[Th]{Definition}
\newtheorem{Pro}[Th]{Proposition}
\newtheorem{Le}[Th]{Lemma}
\newtheorem{Exa}[Th]{Example}
\newtheorem{Rem}[Th]{Remark}
\newtheorem{Fa}[Th]{Fact}
\newtheorem{Que}[Th]{Question}
\newtheorem{Ct}[Th]{}
\newtheorem{Afi}[Th]{Claim}

\newcommand{\baf}{\begin{Afi}\nf{\sl }}
\newcommand{\eaf}{\end{Afi}}

\newcommand{\bdf}{\begin{Df}\nf{\bf }}
\newcommand{\edf}{\end{Df}}
\newcommand{\bte}{\begin{Th}\nf{\bf }}
\newcommand{\ete}{\end{Th}}  
\newcommand{\bco}{\begin{Co}\nf{\bf }}
\newcommand{\eco}{\end{Co}}
\newcommand{\ble}{\begin{Le}\nf{\bf }}
\newcommand{\ele}{\end{Le}}
\newcommand{\bpr}{\begin{Pro}\nf{\bf }}
\newcommand{\epr}{\end{Pro}}
\newcommand{\bex}{\begin{Exa}\nf{\bf } \rm}
\newcommand{\eex}{\end{Exa}}
\newcommand{\bre}{\begin{Rem}\nf{\bf } \rm}
\newcommand{\ere}{\end{Rem}}
\newcommand{\bfa}{\begin{Fa}\nf{\bf } \sl}
\newcommand{\efa}{\end{Fa}}
\newcommand{\bqt}{\begin{Que}\nf{\bf }}
\newcommand{\eqt}{\end{Que}}
\newcommand{\bdm}{{\bf Proof.}}
\newcommand{\edm}{\smallskip}
\newcommand{\qdr}{\hfl $\square$ }
\newcommand{\bct}{\begin{Ct}\nf{\bf } \rm}
\newcommand{\ect}{\end{Ct}}
\newcommand{\bxa}{\begin{Exa}\nf{\bf } \rm}
\newcommand{\exa}{\end{Exa}}

\newcommand{\Ss}{\cS_s}
\newcommand{\Sf}{\cS_f}
\newcommand{\Ls}{\cL_s}
\newcommand{\Lf}{\cL_f}
\newcommand{\Lc}{\cL_f^c}
\newcommand{\Qf}{\cQ_f}
\newcommand{\Qc}{\cQ_f^c}

\begin{document}

\title{Towards a good notion of categories of logics} 

\author{ Caio de Andrade Mendes \thanks{Instituto de Matem\'atica e Estat\'istica, University of S\~ao Paulo, Brazil.\ Emails: caio.mendes@usp.br, hugomar@ime.usp.br }\\
Hugo Luiz Mariano \thanks{Research supported by FAPESP, under the Thematic Project LOGCONS:
Logical consequence, reasoning and computation (number 2010/51038-0).} 
}

\date{March 2016}

\maketitle

\begin{abstract}

We consider (finitary, propositional) logics through the original use
of Category Theory: the study of the `` sociology of mathematical
objects'', aligning us with a recent, and growing, trend of study
logics  through its relations with other logics (e.g. process of 
combinations of logics as fibring \cite{Gab} and possible translation semantics \cite{Car}). So will be objects of
study the classes \emph{of} logics, i.e. categories whose objects
are logical systems (i.e., a signature with a Tarskian consequence relation)  and the morphisms are related to (some concept of)
translations between these systems. 
The present work provides the first steps of a project of considering categories of logical systems satisfying {\em simultaneously} certain natural requirements: it seems that in the literature (\cite{AFLM1},
\cite{AFLM2}, \cite{AFLM3}, \cite{BC}, \cite{BCC1}, \cite{BCC2}, \cite{CG}, \cite{FC})  this is achieved only partially.

\end{abstract}

\section*{Introduction}

\hp We consider (finitary, propositional) logics through the original use
of Category Theory: the study of the `` sociology of mathematical
objects'', aligning us with a recent, and growing, trend of study
logics  through its relations with other logics, e.g. in the process of 
combinations of logics. The phenomenon of combinations of logics (\cite{CC3}), emerged in the mid-1980s, was the main motivation for considering categories of logics. There are two aspects of combination of logics: (i) splitting of logics: a analitical process; (ii) splicing of logics: a synthesis.   
The "Possible-Translations Semantics", introduced in \cite{Car}, is an instance of the splitting process: a given logic system is decomposed into other (simpler) systems, providing, for instance a conservative translation of the logic in analysis into a "product" (or weak product) of simpler or better known logics. The "Fibring" of  logics, introduced originally in the context of modal logics (\cite{Gab}), is "the least logic which extends simultaneously the given logics"; after,  this was recognized as a coproduct construction (\cite{SSC}): this provides an example of synthesis of logics.



In the field of categories of logics there are, of course, two choices that must be done: (i) the choice of objects (how represent a logical system?); (ii) the choice of arrows (what are the relevant notions of morphims between logics?). Here we took very simple and universal choices: a logical system will be a (finitary) signature endowed a Tarskian consequence relation  and the morphisms are related to (some concept of)
"logical translations" between these systems.




The main flow of research on categories of logics, represented by the groups of CLE-Unicamp (Brazil) and IST-Lisboa (Portugal) focus on the determination of the conditions for preservation of metalogical properties under the process of combination of logics (\cite{Con}, \cite{CCCSS}, \cite{CR}, \cite{SRC}, \cite{ZSS}). On the other hand, the "global aspects" of categories of logics, that  ensure for example the abundance or scarcity of constructions,
 seem to have not been adequately studied.


The present work provides the first steps of a project of considering categories of logical systems satisfying {\em simultaneously} certain natural requirements such as:\\
{\bf (i)} If they represent the majority of the usual logical systems; \\
{\bf (ii)} If they have good
categorial properties  (e.g., if they are a complete and/or
cocomplete category, if they are accessible categories (\cite{AR}));\\
{\bf (iii)} If
they allow a natural notion of {\em algebraizable} logical system
(as in the concept of  Blok-Pigozzi algebraizable logic (\cite{BP}) or Czelakowski's proto-algebraizability (\cite{Cze})); \\
{\bf (iv)} If they provide a satisfactory
treatment of the {\em identity problem} of logical systems (when
logics can be considered "the same"? (\cite{Bez}, \cite{CG})).
 
 In the series of articles \cite{AFLM1},
\cite{AFLM2}, \cite{AFLM3}, was considered a  {\bf s}imple (but too {\bf s}trict) notion of morphism of signatures,  where are founded some  categories of logics that satisfy
simultaneously the first three requirements, but  not the item (iv); here we will denote by $\Ss$ and $\Ls$ the category of signatures and of logics therein. 

In the series of papers \cite{BC}, \cite{BCC1}, \cite{BCC2}, \cite{CG}\footnote{We Thank professor Marcelo Coniglio for that reference.}
, \cite{FC} is developed a more {\bf f}lexible notion of morphism of signatures based on {\bf f}ormulas as connectives (our notation for the associated category of signatures will be $\Sf$ and $\Lf$ will denote the associated category of logics),  it encompass itens (i) and (iii) and allows some  treatment of item (iv), but does not satisfy (ii).

In \cite{MM} we provide an approach to overcome  both the   deficiencies of the two series of papers. In the present work we provide some new and more detailed information on the categories of signatures underlying to the categories of logics in the two series of papers above mentioned and also in \cite{MM}: 
We present notions of  categories  of logical
systems (and of of signatures) that do not impose too many constraints and that have not
many categorial failures. We preserve the usual the notion of (finitary, propositional) logic as a pair formed by a (finitary) signature and a  Tarskian consequence relation on the associated set of formulas on denumerable variables,  but we change the notion of (translation) morphism between logics to allow more interesting connections between  logics. The basic idea is to take quotient categories of categories of logics and  translations by a (congruence) relation that identifies two morphims if, for each formula in the domain logic, the associated formulas images by the morphisms in the codomain logic are interdemonstrable, but in fact we work with  reflective subcategory of this quotient category determined by "well-behaved" logics.

We briefly describe the paper. 
Section 1 consider only categories of signatures: in (1.1)  we  recall the basic properties of the categories of signatures  $\Ss$ and $\Sf$ and we add some new information; in (1.2) we compare these two categories of signatures by means of functors \ $\Ss \ \overset{(+)}\to{\underset{(-)}\to\rightleftarrows} \ \Sf$ \ and we prove that they provide an adjoint pair of functors; in (1.3)  we identify the monad (or triple) associated to the described adjunction, we identify some properties of the monad, we prove that $\Sf$ is, precisely, the {\em Kleisli category} of that monad and we extract some consequences.
 Section 2 deals with categories of logics: in (2.1) we describe in details the natural  structure of algebraic lattice on the set of consequence relations over a given signature; in (2.2)  we  recall the basic properties of the categories of logics  $\Ls$ and $\Lf$; in (2.3)  we present  new information on $\Lf$ and we prove that the results on categories of signatures presented in (1.3) and (1.4) "lift" to categories of logics: this constitutes evidence that the defects of $\Lf$ are all inherited from $\Sf$; in (2.4) we introduce  new  categories of logics that solves the  "deficiencies" of the categories of logics presented in the literature.
 We finish the paper in Section 3 with some comments and future perspectives. 

{\em In what follows,  $X = \{x_{0},x_{1},\ldots,x_{n},\ldots\}$
will denote a fixed enumerable set (written in a fixed order).}

\section{Categories of signatures}

\subsection{Known facts about categories of signatures}

\subsubsection{The category $\Ss$}

\hp We will write  $\Ss$ for the category of signatures and {\em strict} morphisms of
signatures presented in \cite{AFLM1}, \cite{AFLM2}, \cite{AFLM3}, and described below. 

The objects of $\Ss$ are signatures. A signature $\Sigma$ is a sequence
of sets $\Sigma = (\Sigma_{n})_{n \in \omega}$ such that $\Sigma_{i} \cap
\Sigma_{j} = \emptyset$ for all $i<j<\omega$ . We write $|\Sigma| 
=  \bigcup_{n \in \omega} \Sigma_{n}$ for the {\em support of $\Sigma$} and we denote by $F(\Sigma)$, the {\em formula algebra of $\Sigma$}, i.e. the set of all
(propositional) formulas built with signature $\Sigma$ over the variables
in $X$.
For all $n \in \N$ let
$F(\Sigma)[n] = \{ \varphi \in F(\Sigma) : var(\varphi) = \{x_{0},x_{1}, \ldots, x_{n-1} \} \}$, where $var(\varphi)$ is the set of all variables that occur in the $\Sigma$-formula $\varphi$.
The notion of complexity $compl(\varphi)$ of the formula $\varphi$ is, as usual, the number of occurences of connectives in $\varphi$.

If $\Sigma, \Sigma'$ are signatures then a  {\em strict} morphism $f : \Sigma \ \lra \
\Sigma'$ is a sequence of functions $f = (f_{n})_{n \in
\omega}$, where $f_{n} : \Sigma_{n} \ \lra \ \Sigma'_{n}$. Composition and identities in $\Ss$ are componentwise.

For each morphism $f : \Sigma \lra \Sigma'$ in $\Ss$ there is a unique  function
$\hat{f} : F(\Sigma) \lra F(\Sigma')$, called the \emph{extension of
$f$}, such that: (i) $\hat{f}(x) = x$, if $x \in X$; (ii) $\hat{f}(c_n(\psi_{0},\ldots, \psi_{n-1}) =
    f_{n}(c_n)(\hat{f}(\psi_{0}),\ldots,\hat{f}(\psi_{n-1}))$, if
    $c_n \in \Sigma_{n}$. Then, by induction on the complexity of formulas:\\
{\bf (0)} \ $compl(\hat{f}(\theta)) = compl(\theta))$, for all $\theta \in F(\Sigma)$.\\
{\bf (1)} \
If $var(\theta) \ \sub \ \{x_{i_0}, \ldots, x_{i_{n-1}}\}$, then  $\hat{f}(\theta(\vec{x}) [\vec{x} \mid \vec{\psi}]) = (\hat{f}(\theta(\vec{x})) [\vec{x} \mid \hat{f}(\vec{\psi})]$. Moreover $var(\hat{f}(\theta)) = var(\theta)$ and then $\hat{f}$ restricts to maps $ \hat{f}\rest_n : F(\Sigma)[n] \ \lra \ F(\Sigma')[n]$, \ $n \in \N$. \\
{\bf (2)} \ The extension to the formula
algebra of a composition is the extension's composition. The extension of an identity is the identity function
on the formula algebra.

Remark that $\Ss$ is equivalent to the functor category
$\mathbf{Set}^\mathbb{N}$, where $\mathbb{N}$ is the discrete category with
object class $\N$, then $\cS$ has all small limits and colimits and
they are componentwise.
Moreover, the category $\Ss$ is a finitely locally presentable category, i.e., $\Ss$ is a finitely 
accessible category that is cocomplete and/or complete(\cite{AR}).
The finitely presentable signatures are precisely the signatures of finite support.

\nhp {\bf (Sub)}  \  For any \emph{substitution} function $\sigma : X \lra F(\Sigma)$,
    there is unique \emph{extension}
    $\widetilde{\sigma}: F(\Sigma) \lra F(\Sigma)$
    such that
    $\widetilde{\sigma}$ is an ``homomorphism'':
    $\widetilde{\sigma}(x) = \sigma(x)$, for all $x \in X$ and
    $\widetilde{\sigma}(c_{n}(\psi_{0}, \ldots, \psi_{n-1}) =
    c_{n}(\widetilde{\sigma}(\psi_{0})), \ldots,
    \widetilde{\sigma}(\psi_{n-1}))$, for all $c_{n} \in \Sigma_{n}$,
    $n \in \omega$;
    it follows that for any $\theta(x_{0},\ldots, x_{n-1}) \in
    F(\Sigma)$ \ $\td{\sigma}(\theta(x_{0},\ldots, x_{n-1})) =
    \theta(\sigma(x_{0}),\ldots, \sigma(x_{n-1}))$.
    The \emph{identity substitution} induces the identity homomorphism
    on the formula algebra; the {\em composition substitution} of the
    substitutions $\sigma', \sigma : X \lra F(\Sigma)$ is the
    substitution $\sigma'' : X \lra F(\Sigma)$ ,  $\sigma'' =
    \sigma'\star\sigma := \td{\sigma'} \circ \sigma$ and
    $\td{\sigma''} = \td{\sigma'\star\sigma} = \td{\sigma'} \circ
    \td{\sigma}$.

\nhp {\bf (3)}  \ Let $f : \Sigma \lra \Sigma'$ be a $\Ss$-morphism.  Then for each
    substitution $\sigma : X \lra F(\Sigma)$  there is a 
    substitution $\sigma' : X \lra F(\Sigma')$ such that
    $\widetilde{\sigma'}\circ \hat{f} = \hat{f} \circ
    \widetilde{\sigma}$.

\subsubsection{The category $\Sf$}

\hp We will write  $\Sf$ for the category of signatures and {\em flexible} morphisms of
signatures presented in the series of papers \cite{BC}, \cite{BCC1}, \cite{BCC2}, \cite{CG}, \cite{FC}  and described below.

We introduce the following notations:\\
If $\Sigma = (\Sigma_n)_{n \in \N}$ is a signature, then write $T(\Sigma) := (F(\Sigma)[n])_{n \in \N}$; clearly $T(\Sigma)$ satisfies the "disjunction condition", then it is a signature too.\\
We have the inverse bijections (just notations): \\
 $h \in \Sf(\Sigma, \Sigma') \ \lrsa \ h^\sharp \in \Ss(\Sigma, T(\Sigma'))$ \hfl \hfl
$f \in \Ss(\Sigma, T(\Sigma')) \ \lrsa \ f^\flat \in \Sf(\Sigma, \Sigma')$.\\
For each signature $\Sigma$ and $n \in \N$, let the function:\\
 $(j_{\Sigma})_n : \Sigma_{n} \ \lra \ F(\Sigma)[n]$ \hem  : \hem $c_n$ \ $\mapsto$ \ $c_n(x_0, \ldots, x_{n-1})$.

For each morphism $f : \Sigma \lra \Sigma'$ in $\Sf$ there is a unique  function
$\check{f} : F(\Sigma) \lra F(\Sigma')$, called the \emph{extension of
$f$}, such that: (i) $\check{f}(x) = x$, if $x \in X$; (ii) $\check{f}(c_n(\psi_{0},\ldots, \psi_{n-1}) =
    (f_{n}(c_n)(x_0, \ldots, x_{n-1}))[x_0 \mid \check{f}(\psi_{0}),\ldots,x_{n-1} \mid \check{f}(\psi_{n-1}))$, if
    $c_n \in \Sigma_{n}$.

The notion of extension of $\Sf$-morphism to formula algebras shares many properties with notion of extension of $\Ss$-morphism to formula algebras: e.g., the properties {\bf (1)}, {\bf (2)}, {\bf (3)}.

The composition in $\Sf$ is given by $(f' \bullet f)^\sharp := (\check{f'}\!\!\rest_{n} \circ (f^\sharp)_{n})_{n \in \N}$. The identity $id_{\Sigma}$ in $\Sf$ is given by $id_{\Sigma}^{\sharp} = ((j_{\Sigma})_{n})_{n \in \N}$. \\



Remark that the "information encoded" by the of extension of $\Sf$-morphism is enough to determine that morphism.
More precisely, given $g, f \in \Sf(\Sigma,\Sigma')$,  note that:\\
$\ast$\ $(f^\sharp)_{n} = \check{f}\!\!\rest_{n} \circ (j_\Sigma)_{n}$, $n \in \N$;\\
$\ast$ \ $\check{f} = \check{g} \ \im \ f^\sharp = (\check{f}\!\!\rest_{n})_{n \in \N} \circ j_\Sigma = (\check{g}\!\!\rest_{n})_{n \in \N} \circ j_{\Sigma} = g^\sharp \ \im\ f = g$.

 For the reader's convenience we add here the proof that $\Sf$ is a category:\\
  $\ast$ \ \und{identity}: $f \bullet id_{\Sigma} = f = id_{\Sigma'} \bullet f$:\\
 $(f \bullet id_{\Sigma})_n^\sharp \ =\ \check{f}\!\!\rest_{n} \circ (id_{\Sigma})_n^\sharp \ =\ \check{f}\!\!\rest_{n} \circ (j_{\Sigma})_n = (f^\sharp)_n$; \\
  $(id_{\Sigma'} \bullet f)_n^\sharp\ =\  \check{id_{\Sigma'}}\!\!\rest_{n}\circ (f^\sharp)_n \ =\ id_{F(\Sigma')[n]} \circ (f^\sharp)_n \ = \ (f^\sharp)_n$.\\
 $\ast$ \ \und{associativity}: $(f''\bullet f')\bullet f = f''\bullet (f'\bullet f)$:\\
 $({(f''\bullet f')\bullet f})_n^\sharp  \ =\ {\check{(f'' \bullet f')}\!\!\rest_{n}}\circ (f^{\sharp})_n \ =\ (\check{f''} \circ \check{f'})\!\!\rest_{n}\circ 
  (\check{f}\!\!\rest_{n} \circ (j_\Sigma)_{n})$ \ = \\
   $\check{f''}\!\!\rest_{n} \circ (\check{f'}\!\!\rest_{n}\circ  
 \check{f}\!\!\rest_{n}) \circ (j_\Sigma)_{n} \ =\ 
   \check{f''}\!\!\rest_{n} \circ (\check{f'}\circ  
  \check{f})\!\!\rest_{n}) \circ (j_\Sigma)_{n} \ =\
  \check{f''}\!\!\rest_{n} \circ (\check{f'\bullet f})\!\!\rest_{n} \circ (j_\Sigma)_{n} \ =\   
  \check{f''}\!\!\rest_{n} \circ ((f' \bullet f)^\sharp)_n \ =\ 
 (f''\bullet (f'\bullet f))_n^\sharp$.

The notion of extension of $\Sf$-morphism to formula algebras shares many properties with notion of extension of $\Ss$-morphism to formula algebras, however: \\
 {\bf (0)}  \ 
 If $f \in \Sf(\Sigma, \Sigma')$, then are equivalent: 
 \\ $\ast$ \ $compl(\check{f}(\theta))  \geq compl(\theta)$, any $\theta \in F(\Sigma)$;\\
 $\ast$ \ $f(c_1) \neq x_0$ , all $c_1 \in \Sigma_1$.

Now we add  some information easily established:

\bdf \label{regmor-df} A $\Sf$-morphism $f : \Sigma \lra \Sigma'$ is {\em regular} if \ $compl(\check{f}(\theta))  \geq compl(\theta))$, any $\theta \in F(\Sigma)$. 
\qdr\edf

\bpr \label{reg-pr} (a) \ If $f \in \Sf(\Sigma, \Sigma')$, then: $f$ is regular \ iff \ $f(c_1) \neq x_0$ , all $c_1 \in \Sigma_1$.\\
(b) \ The "empty" signature is the unique initial object of $\Sf$ (as in $\Ss$).\\
(c) \ A (non full) subcategory of $\Sf$ with the same objects has  a {\em strict } initial object iff all the morphisms are regular.\\
(d) The mapping \ $f \in \Ss(\Sigma, \Sigma') \ \mapsto (j_{\Sigma'} \circ f)^\flat\in \Sf(\Sigma, \Sigma')$ is a  {\em (natural) bijection}   $\Ss(\Sigma, \Sigma') \ \overset{\cong}\to\lra \  \{h \in\Sf(\Sigma, \Sigma') : compl(\check{h}(\theta))  = compl(\theta)$, for any $\theta \in F(\Sigma)\}$.
\epr
\bdm The non trivial implication in item (a) follows from induction on complexity of formulas.
\edm\qdr

\bpr \label{terminal-pr} (a) $\Sf$ has  weak terminal objects . More precisely, a signature $\Sigma$ is an weak terminal object iff  $\Sigma'_0 \neq \empty$ and exists $k \geqs 2$ such that $\Sigma'_k \neq\empty$. \\
(b) $\Sf$ does not have terminal object. (Example: let $\Sigma'$ be a weak terminal object and take a signature $\Sigma$ with only one conective and it is binary: as $F(\Sigma')[2]$ is {\em infinity}, there are many $\Sf$-morphisms from $\Sigma$ into $\Sigma'$.)
\qdr
\epr

\bre \label{weakprodequal-re}\\
(a)  It is easy to see that $\Sf$ has weak products: a weak product of a (small) family of signatures can be given by taking the product signature in the {\em strict category} $\Ss$ and the corresponding $\Ss$-projections, transformed into $\Sf$-morphisms (see the next subsection).\\
(b)  As $\Sf$ has initial object, any family of paralel arrows has an weak equalizer.\\  
(c)$\Sf$ has  weak terminal object but do not have terminal object: a necessary and sufficient condition to $\Sigma$  be a weak terminal 
is $\Sigma_0 \neq \empty, exists k \geqs 2 Sigma_k \neq\empty$
(d) $\Sf$ only has "trivial" (i.e. in $\Ss$) coequalizers, idempotents, isomorphisms,
\qdr\ere

\subsection{The fundamental adjunction}

\bpr\label{coneccatsig-pr} {\em Conecting categories of signatures:}\\
(a) We have the (faithful) functors: \\ 
$(+) : \Ss \ \lra \ \Sf$ \hem : \hem $(\Sigma \overset{f}\to\lra\ \Sigma')$ \ $\mapsto$ \ $(\Sigma \overset{(j_{\Sigma'} \circ f)^\flat}\to\lra\ \Sigma')$;\\
$(-) : \Sf \ \lra \ \Ss$ \hem : \hem $(\Sigma \overset{h}\to\lra\ \Sigma')$ \ $\mapsto$ \ $((F(\Sigma)[n])_{n \in \N} \overset{(\check{h}\rest_n)_{n\in\N}}\to\lra\ (F(\Sigma' )[n])_{n \in \N})$.\\
(b) For each $f \in \Ss(\Sigma, \Sigma')$, we have $\check{(f^+)} = \hat{f} \in Set(F(\Sigma), F(\Sigma'))$.\\
(c) We have the {\em natural transformations}: \\
 $\eta : Id_{\Ss} \ \lra \ {(-) \circ (+)}$ \ : \ $(\eta_{\Sigma})_{n} := (j_{\Sigma})_n$ \hfl \hfl
$\varepsilon : (+) \circ (-) \ \lra \ Id_{\Sf}$ \ : \ $(\varepsilon_{\Sigma})^\sharp_{n} := id_{F(\Sigma)[n]}$ \\
 and we write
$\mu = (-)\varepsilon(+)$. 
\qdr\epr

\bte \label{adjun-te} The (faithful) functor $(+)$ is a {\em left adjoint} of the (faithful) functor $(-)$: $\eta$ and $\varepsilon$ are, respectively, the unit and the counit of the adjunction.
\qdr\ete


\bco \label{} (a)\ The functor $(+)$ preserves colimits and the functor $(-)$ preserves limits.\\
\qdr\eco

\bco\ \label{Sfcolim-co} The category $\Sf$ has colimits for any (small) diagram "in $\Ss$", i.e., given $\cI$ a small category and a diagram $D : \cI \ \lra \ \Ss$, the category $\Sf$ has a colimit for the diagram $(+) \circ D : \cI \ \lra \ \Sf$. In particular, $\Sf$ has all (small) coproducts and all (small) pushouts "based in $\Ss$".
\qdr \eco

\bpr \label{} Let $h \in \Sf(\Sigma, \Sigma')$:\\
(a) \ If $h^{-}$ is a $\Ss$-epimorphism, then $f$ is $\Sf$-epimorphism.\\
(b) \ $h$ is a $\Sf$-monomorphism if and only if  $h^{-}$ is a $\Ss$-monomorphism.\\
(c) If $h$ is an $\Sf$-isomorphism, then  ``$h \in \Ss$'', i.e. there is a (unique) $\Ss$-(iso)morphism $f$ such that $h = f^-$; in particular, $h$ is regular. \\
(d) If $h$ is a $\Sf$-section, then $h$ is regular and if $g\bullet h =id$ for some $\Sf$-morphism $g$ that is regular over the "image signature of $h$" (i.e. the signature whose conectives effectively occur in  the image of some $h_n $, $n \in \N$), then 
``$h \in \Ss$''.
\qdr\epr

\subsection{The monad and its properties} 

\hp We have a (endo)functor $T : \Ss \ \lra\ \Ss$    \hem $(\Sigma \overset{f}\to\lra\ \Sigma')$ \ $\overset{T}\to\mapsto$ \
$((F(\Sigma)[n])_{n \in \N} \overset{(\hat{f}\!\rest_n)_{n\in\N}}\to\lra\ (F(\Sigma' )[n])_{n \in \N})$. Clearly, $T = (-) \circ (+)$ and it is a faithful functor. Let  $\cT = (T, \eta, \mu)$ be the monad (or triple)  associated  to the  adjunction $(\eta, \varepsilon) : \Ss \ \overset{(+)}\to{\underset{(-)}\to\rlas} \ \Sf$.

\bpr \label{Tiso-pr} The functor $T$ reflects  isomorphisms (respectively: monomorphisms, epimorphisms).  
\epr
\bdm First remark that, for each signature $\Sigma$ and $n \in\N$, $(\eta_\Sigma)_n : \Sigma_n \rat F(\Sigma)[n]$ stablish a bijection between $\Sigma_{n}$ and $\{ \theta \in F(\Sigma)_{n}: compl(\theta)=1\}$. Now let $f : \Sigma \ \lra \Sigma'$ a $\Ss$-morphism such that $T(f)$ is a $\Ss$-isomorphism (respectively: a $\Ss$-monomorphism, a $\Ss$-epimorphism). Then, for each $n \in \N$, $\hat{f}\!\!\rest_{n} :  F(\Sigma)[n] \ \lra \ F(\Sigma')[n]$ is a bijection (respectively: a injection, a surjection) and, as $compl(\hat{f}(\theta)) = compl(\theta)$ for each $\theta \in F(\Sigma)$, $\hat{f}\!\!\rest_{n}$ restricts to a bijection (respectively: a injection, a surjection) between $\{ \theta \in F(\Sigma)[n]: compl(\theta)=1\}$ and $\{ \theta' \in F(\Sigma')[n]: compl(\theta')=1\}$. Finally, as $\eta_{\Sigma'} \circ f = T(f) \circ \eta_{\Sigma}$, we conclude that $f_{n} : \Sigma_{n} \ \lra \ \Sigma'_{n}$ is a bijection (respectively: a injection, a surjection), for each $n \in N$, as we need.
\qdr\edm

                                         

\bpr \label{Tdircolim-pr} The functor $T$  preserves directed colimits (i.e., colimits of diagrams over upward directed posets). More explicitly, let  $(I,\leqs)$  be an upward directed poset and $D : (I,\leqs) \ \lra \ \Ss$ : $i \mapsto \Sigma_i$ be  a diagram in $\Ss$; $(\Sigma, (\Sigma_i \overset{\alpha_i}\to\ra \Sigma)_{i \in I})$  denotes the colimit of $D$ in $\Ss$; $(\Sigma', (T(\Sigma_i) \overset{\alpha'_i}\to\ra \Sigma')_{i \in I})$ be the colimit of $T\circ D$ in $\Ss$; $(S, (F(\Sigma_i) \overset{\beta_i}\to\ra S)_{i \in I})$ denotes the colimit of $\hat{(\ )} \circ D$ in the category $Set$, then:\\
(a) \ The canonical function  $S \ra F(\Sigma)$, denoted $k : colim_{i \in I} F(\Sigma_i) \ \lra \  F(colim_{i \in I} \Sigma_i)$,  i.e. the unique function such that $k \circ \beta_i = \hat{\alpha}_i$, $i \in I$, is a {\em bijection}.\\
(b)  \ The canonical $\Ss$-morphism    $can  : colim_{i \in I} T(\Sigma_i) \ \lra \  T(colim_{i \in I} \Sigma_i)$, i.e.  the unique $\Ss$-morphism such that $can \circ \alpha'_i = T(\alpha_i)$, $i \in I$, is a \und{$\Ss$-isomorphism}. It is given by sequence of {\em bijections} \ $can_n  : colim_{i \in I} (F(\Sigma_i)[n]) \ \lra$ \
$  F(colim_{i \in I} \Sigma_i)[n]$, $n \in \N$, obtained from the  "restrictions" of the canonical bijection  $k$ just above. 
\epr
\bdm ({\em Sketch}) For a proof of item (a) we apply a "global reasoning": we consider formula algebras and apply  induction on complexity of formulas. For (b): we extract  "local" information from the result is (a), i.e., we  consider convenient "restrictions" to the subsets  $F(\Sigma)[n]$, $n \in \N$.
\qdr\edm

The same technique of proof in the Proposition above gives us the Theorem below:

\bte \label{Kleisli-te} Let  $\cT = (T, \eta, \mu)$ be the monad associated to the adjunction $ (\eta, \varepsilon) : \Ss \ \overset{(+)}\to{\underset{(-)}\to\rlas} \ \Sf$ (i.e., $\mu = (-)\varepsilon(+)$)  is such that $Kleisli(\cT) = \Sf$. Moreover, the functors $(+)$ and $(-)$ are precisely the canonical functors associated to the adjunction  of the Kleisli category of a monad. More explicitly: given $(\Sigma \ \overset{f}\to\lra \ \Sigma' \ \overset{f'}\to\lra \ \Sigma'')$ in $\Sf$, then $f' \bullet f = (\mu_{\Sigma''} \circ T(f'^{\sharp}) \circ f^\sharp)^\flat$ ,i.e.,  we have in $\Ss$:\\
 \hfl $(\Sigma \ \overset{f^\sharp}\to\lra \ T(\Sigma') \ \overset{(\check{f'}\rest_{n})_{n\in \N}}\to\lra \ T(\Sigma''))$ \hem = \hem $(\Sigma \ \overset{f^\sharp}\to\lra \ T(\Sigma') \ \overset{T(f'^{\sharp})}\to\lra \ T \circ T(\Sigma'') \ \overset{\mu_{\Sigma''}}\to\lra \ T(\Sigma''))$.\hfl
 \qdr\ete


\section{Categories of Logics}

\subsection{The lattice of logics above a signature}

\hp A logic is an ordered pair $l = (\Sigma,
\vdash)$ where $\Sigma$ is an object of $\Ss$ and $\vdash$
codifies the (Tarskian) ``consequence operator'' on $F(\Sigma)$ \ : \ $\vdash$ is
a binary relation, a subset of $Parts(F(\Sigma)) \times F(\Sigma)$, such
that $Cons(\Gamma) = \{ \varphi \in F(\Sigma) : \Gamma \vdash
\varphi\}$, for all $\Gamma \sub F(\Sigma)$, gives a structural finitary
closure operator on $F(\Sigma)$: \\
{\bf (a)} {\em inflationary}: $\Gamma \sub Cons(\Gamma)$; \\
{\bf (b)} {\em increasing}: $\Gamma_{0} \sub \Gamma_{1} \im Cons(\Gamma_{0})
    \sub Cons(\Gamma_{1})$;\\
{\bf (c)} {\em idempotent}: $Cons(Cons(\Gamma)) \sub Cons(\Gamma)$;\\
{\bf (d)} {\em finitary}: $Cons(\Gamma) = \bcup \{ Cons(\Gamma') : \Gamma'
    \sub_{fin} \Gamma\}$; \\
{\bf (e)} {\em structural}: $\widetilde{\sigma}(Cons(\Gamma)) \sub
    Cons(\widetilde{\sigma}(\Gamma))$, for each substitution $\sigma :
    X \to F(\Sigma)$.
    
The set  of all consequence relations on a signature $\Sigma$, denoted by
$Cons_{\Sigma}$, is endowed with the partial order: \ $\vdash_0 \leqs \vdash_1$ \ iff \ for each $\Gamma \in Parts(F(\Sigma))$, $\overline{\Gamma}^0 \sub \overline{\Gamma}^1$.

\bfa\ For each signature $\Sigma$, the poset  $(Cons_{\Sigma}, \leqs)$ is a complete lattice. It is in fact an algebraic lattice
where the compact elements are the ``finitely generated logics'', the
logics over $\Sigma$ given by a finite set of axioms and a finite set of
(finitary) inference rules.
\efa
\bdm (sketch)

\und{Infs}:
Consider $I$ a set and $D = \{l^{i}=(\Sigma,\vdash_{i})\}_{ i \in I}$ a family
 of logics over the signature $\Sigma$. Now, for each $\Gamma \cup\{\psi\} \ \sub
 \ F(\Sigma)$, define that $\Gamma \vdash \psi \Lra$ there is
 $\Gamma'\sub_{fin}\Gamma$ such that $(\forall i \in I)(\Gamma' \vdash_{i}
 \varphi)$, then $(\Sigma,\vdash)$ is a logic and $l = (\Sigma,\vdash)$ is the
 infimum of the family $D$ in $\cL$, thus $(Cons_\Sigma, \leqs)$ is a {\em complete lattice}.

 \und{Generated consequence relation}:
 As the set of consequence operators (or consequence relations) on a signature
 $\Sigma$ is a complete lattice, there exists the logic generated by any function
 $W : Parts(F(\Sigma)) \lra Parts(F(\Sigma))$: it is enough to take the infimum of the
 family of all consequence relations on $\Sigma$ that are upper bounds of the
 ``proto-consequence relation'' $W\!\!\rest : Parts_{fin}(F(\Sigma)) \lra Parts(F(\Sigma))$ associated with  $W$. A more explicity description is given by the usual notion of "proof"  based on hypothesis, axioms and inference rules.


 \und{Directed sups}:
 Consider $I$ a set and $D = \{l^{i}=(\Sigma,\vdash_{i})\}_{ i \in I}$ an
 upward directed family of logics over the signature $\Sigma$, that is, for each
 $i,j \in I$ there is a $k \in I$ such that $id_{\Sigma} \in \cL(l^{i},l^{k})$ ,
 $id_{\Sigma} \in \cL(l^{j},l^{k})$. Now, for each $\Gamma \cup\{\psi\} \sub
 F(\Sigma)$, define that $\Gamma \vdash \psi \Lra$ there is
 $\Gamma'\sub_{fin}\Gamma$ and there is an $i \in I$ such that $\Gamma'
 \vdash_{i} \varphi$, then $(\Sigma,\vdash)$ is a logic and $l = (\Sigma,\vdash)$
 is the supremum of the family $D$ in $\cL$.  

 \und{Sups}:
 As usual, the supremum of a family of logics can be obtained taking the infimum
 of the set of upper bounds of that family of logics. A more objective
 characterization of suprema can be given but we postpone that because this can
 be easily described by more general results below (see
 Proposition~\ref{Llim-pr}).

 \und{Compact consequence relations}:
 A consequence relation $\vdash'$ over $\Sigma$ is compact if for
 each set $I$, each $D = \{l^{i}=(\Sigma,\vdash_{i})\}_{ i \in I}$ a
 upward directed family of logics over the signature $\Sigma$, if
 $\vdash' \leq \bigvee_{i \in I} \vdash_{i}$ then there is an $i \in
 I$ such that $\vdash' \leq \vdash_{i}$.  It follows easily that this
 condition is equivalent to the ``stronger'' condition: for each set
 $J$, each $D = \{ l^{j}=(\Sigma,\vdash_{j})\}_{ j \in J}$ a family
 of logics over the signature $\Sigma$, if $\vdash' \leq \bigvee_{j
 \in J} \vdash_{j}$ then there is a finite subset $J' \sub J$ such
 that $\vdash'\ \leq\ \bigvee_{j \in J'} \vdash_{j}$.\footnote{Just
 observe that any sup of a family coincides with a sup of a directed
 family: for each set $J$ take $I = P_{fin}(J)$ then, for each
 $J'\sub_{fin} J$, define $\vdash_{J'} = \bigvee_{j \in J'}
 \vdash_{j}$ \ldots .} A consequence relation on $\Sigma$ is compact
 if and only if it is a finitely generated consequence relation on
 $\Sigma$. Any consequence relation on $\Sigma$ is the directed
 supremum of its compact (sub)consequence relations on $\Sigma$. 
 \qdr\edm

\subsection{Known facts about categories  of logics}


\hp The category $\Ls$ is the category of propositional logics and {\em strict} translations
as morphisms.  This is a category ``built above'' the category $\Ss$, that
is, there is an obvious forgetful functor $U_s : \Ls \ \lra \ \Ss$. The categorial properties of $\Ls$ are detailed in \cite{AFLM3}.

The objects of $\Ls$ are logics $l= (\Sigma, \vdash)$ as described in subsection 2.1.

If $l=(\Sigma,\vdash), l'=(\Sigma',\vdash')$ are logics then a
\emph{strict translation morphism} $f : l \lra l'$ in $\Ls$ is a {\em strict} signature morphism
$f : \Sigma \lra \Sigma'$ in $\Ss$ such that ``preserves the consequence relation'', that
is, for all $\Gamma \cup\{\psi\} \sub F(\Sigma)$, if $\Gamma \vdash \psi$
then $\widehat{f}[\Gamma] \vdash' \widehat{f}(\psi)$.  Composition and
identities are similar to $\Ss$.

$\Ls$ has  natural notions of direct and inverse image logics under a $\Ss$-morphism and they have good properties:

\bct\  {\em Direct image and inverse image:}\\
Let $f : \Sigma \lra \Sigma'$ be a
$\Ss$-morphism:\\
 \textbf{Inverse image:} if $l'=(\Sigma',\vdash') \in Obj(\cL)$
    then for all $\Gamma \cup\{\psi\} \sub F(\Sigma)$ define
    $\Gamma \vdash_{f^{\star}(\vdash')}\psi$ iff $\widehat{f}[\Gamma]
    \vdash' \widehat{f}(\psi)$;\\
 \textbf{Direct image:} if $l=(\Sigma,\vdash) \in Obj(\cL)$ then
    for all $\Gamma' \cup\{\psi'\} \sub F(\Sigma')$ define $\Gamma'
    \vdash_{f_{\star}(\vdash)} \psi'$ iff there is a finite sequence
    of $\Sigma'$-formulas $(\phi'_{0},\ldots,\phi'_{t})$ such that:\\
  $\bullet$ \ $\phi'_{t} = \psi'$;\\
    $\bullet$ \ for all $p \leq t$ at least one of the alternatives below
        occurs:\\
       \hp $\ast$\ ``$\phi'_{p}$ is a hypothesis'': $\phi'_{p} \in
            \Gamma'$;\\
          \hp $\ast$\  ``$\phi'_{p}$ is an instance of an $l$-axiom'': there
            is a $\theta_{p} \in F(\Sigma)$ such that $\vdash
            \theta_{p}$ and there is a substitution $\sigma' : X
            \lra F(\Sigma')$ such that
            $\widetilde{\sigma}'(\widehat{f}(\theta_{p}))=
            \phi'_{p}$;\\
          \hp $\ast$\  ``$\phi'_{p}$ is a direct consequence of an instance
            of $l$-inference rule applied over previous members in
            the sequence'': there is a $\Delta_{p} \cup
            \{\theta_{p}\} \sub_{fin} F(\Sigma)$ such that
            $\Delta_{p}\vdash \theta_{p}$ and there is a
            substitution $\sigma' : X \lra F(\Sigma')$ such that
            $\widetilde{\sigma}'(\widehat{f}(\theta_{p}))=
            \phi'_{j}$ and
            $\widetilde{\sigma}'[\widehat{f}[\Delta_{p}]] \sub
            \{\phi'_{0}, \ldots, \phi'_{j-1}\}$.
   \qdr \ect

\bfa \label{image-fa}%
 Let $f : \Sigma \lra
\Sigma'$ be a $\Ss$-morphism and let $l=(\Sigma,\vdash),
l'=(\Sigma',\vdash')$ be logics $l,l' \in Obj(\Ls)$. Then
$(i)^{\star}$ and $(i)_{\star}$ hold, and $(ii)^{\star}, (ii)_{m}$
and $(ii)_{\star}$ are equivalent:\\
$(i)^{\star}$ \ if $l'=(\Sigma',\vdash') \in Obj(\Ls)$ then
    $f^{\star}(l') = (\Sigma, \vdash_{f^{\star}(\vdash')}) \in Obj(\Ls)$;\\
$(i)_{\star}$ \ if $l=(\Sigma,\vdash) \in Obj(\Ls)$ then
    $f_{\star}(l) = (\Sigma', \vdash_{f_{\star}(\vdash)}) \in
    Obj(\Ls)$.\\
$(ii)^{\star}$ \ $\vdash\ \leq\ f^{\star}(\vdash')$;\\
$(ii)_{m}$ \ $f : (\Sigma,\vdash) \lra (\Sigma',\vdash')$ is $\Ls$-morphism;\\
$(ii)_{\star}$\ $f_{\star}(\vdash)\ \leq'\ \vdash'$.
\qdr\efa

\bre\ \label{adjunts-re}%
It follows easily from the facts above that the forgetful functor
$U_s : \Ls \lra \Ss$  : $((\Sigma, \vdash)$ $\overset{f}\to\lra (\Sigma',\vdash'))
\mapsto (\Sigma \overset{f}\to\lra \Sigma')$ has left and right adjoint
functors: the left adjoint  $\bot_s : \Ss \lra \Ls$ and the right adjoint $\top_s :
\Ss \lra \Ls$ take a signature $\Sigma$ to, respectively, $\bot_s(\Sigma) =
(\Sigma, \vdash_{min})$ (the first element of $Cons_{\Sigma}$) and
$\top_s(\Sigma) = (\Sigma, \vdash_{max})$ (the last element of $Cons_{\Sigma}$). Moreover, $U_s \circ \bot_s = Id_{\Ss} = U_s \circ \top_s$ and $U_s$ preserves all limits and colimits that exists in $\Ss$.
\qdr\ere

\bfa\ \label{Llim-pr}%
The category $\Ls$ is complete and cocomplete and the forgetful functor
$U_s : \Ls \lra \Ss$ "lifts" all small limits and colimits.
\qdr \efa

\bct\ \label{Lslimcolim-co} We  describe explicitly limits, directed colimits and (general) colimits in $\Ls$:

\und{Limits}:
Let $\cI$ be a small category and  $D : \cI \lra \cL$,
$((\Sigma^{i},\vdash_{i}) \overset{f^{h}}\to\lra (\Sigma^{j},\vdash_{j}))_{(h:
i \to j)\in \cI}$ a diagram, and take $(\Sigma,(\pi^{i})_{i \in Obj(\cI)})$
the limit of the underlying diagram $(\cI \overset{D}\to\lra \cS
\overset{U}\to\lra \cL)$. For all $\Gamma \cup\{\psi\} \sub F(\Sigma)$, define
that $\Gamma \vdash \psi \Lra$ there is ${\Gamma^{-}}\sub_{fin}\Gamma$ such
that for all $i \in Obj(\cI)$ \ $\widehat{\pi}^{i}[{\Gamma^{-}}]
\vdash_{i}\widehat{\pi}^{i}(\psi)$\footnote{This definition also works for
the terminal logic $l=(\Sigma,\vdash)$ where $\Sigma$ is the terminal
signature ($card(\Sigma_{n})=1, \all n \in \omega$), and for all $\Gamma
\cup\{\psi\} \sub F(\Sigma)$, $\Gamma \vdash \psi$.}, then
$l=(\Sigma,\vdash)$ is a logic and $(l,(\pi^{i})_{i \in Obj(\cI)})$ is the
limit of $D$ in $\cL$.

\und{Directed colimits}:
Let $(I,\leq)$ be a directed ordered set and $D : (I,\leq) \lra \cL$,\\
$((\Sigma^{i},\vdash_{i}) \overset{f^{ij}}\to\lra (\Sigma^{j},\vdash_{j}))_{(i
\leq j)\in I}$ be a diagram. Take $(\Sigma,(\gamma^{i})_{i \in I})$ the
colimit of the underlying diagram $(\cI \overset{D}\to\lra \cS \overset{U}\to\lra
\cL)$. Now, for all $\Gamma \cup\{\psi\} \ \sub \ F(\Sigma)$, define that
$\Gamma \vdash \psi \Lra$ there is ${\Gamma^{-}}\sub_{fin}\Gamma$ and there
is an $i \in I$ such that ${\Gamma^{-}} \cup\{\psi\} \sub
\widehat{\gamma}^{i}[F(\Sigma^{i})]$ and there is ${\Gamma^{-}}^{i}
\cup\{\psi^{i}\} \sub_{fin} F(\Sigma^{i})$ such that
$\widehat{\gamma}^{i}[{\Gamma^{-}}^{i}]={\Gamma^{-}}$ ,
$\widehat{\gamma}^{i}(\psi^{i})= \psi$ and \ ${\Gamma^{-}}^{i}
\vdash_{i}\psi^{i}$. Then $l=(\Sigma,\vdash)$ is a logic and
$(l,(\gamma^{i})_{i \in I})$ is the colimit of $D$ in $\cL$.

\und{Colimits}:
Let $\cI$ be a small category and  $D : \cI \lra \cL$,
$((\Sigma^{i},\vdash_{i}) \overset{f^{h}}\to\lra (\Sigma^{j},\vdash_{j}))_{(h:
i \to j)\in \cI}$ be a diagram, and take $(\Sigma,(\gamma^{i})_{i \in
Obj(\cI)})$ the colimit of the underlying diagram $(\cI \overset{D}\to\lra \cS
\overset{U}\to\lra \cL)$. Now, for all $\Gamma \cup\{\psi\} \sub F(\Sigma)$,
define that $\Gamma \vdash \psi \Lra$ there is a finite sequence of
$\Sigma$-formulas $(\phi_{0}, \ldots, \phi_{t})$, where $\phi_{t} = \psi$
and for all $p\leq t$ one of these alternative occurs:\\
  $\bullet$ \``$\phi_{p}$ is an hypothesis'': $\phi_{p} \in \Gamma$;\\
   $\bullet$ \``$\phi_{p}$ is an axiom'': there are $i \in Obj(\cI)$,
    $\theta^{i} \in F(\Sigma^{i})$, $\sigma: X \lra F(\Sigma)$ such
    that $\vdash_{i} \theta^{i}$ and $\phi_{p} =
    \widetilde{\sigma}(\widehat{\gamma}^{i}(\theta^{i}))$;\\
   $\bullet$ \ ``$\phi_{p}$ is a consequence of a inference rule'': there are $i
    \in Obj(\cI)$, $\Delta^{i}\cup\{\theta^{i}\} \sub_{fin}
    F(\Sigma^{i})$, $\sigma: X \lra F(\Sigma)$ such that
    $\Delta^{i}\vdash_{i} \theta^{i}$ and
    $\widetilde{\sigma}[\widehat{\gamma}^{i}[\Delta^{i}]] \sub
    \{\phi_{0},\ldots,\phi_{p-1}\}$, $\phi_{p} =
    \widetilde{\sigma}(\widehat{\gamma}^{i}(\theta^{i}))$;\\
Then $l=(\Sigma,\vdash)$ is a logic and $(l,(\gamma^{i})_{i \in Obj(\cI)})$
is the colimit of $D$ in $\cL$.
\qdr\ect

\bfa\ The category $\Ls$ is a finitely locally presentable category, i.e., $\Ls$ is a finitely 
accessible category that is cocomplete and/or complete (\cite{AR}). The finitely presentable objects in $\Ls$ are precisely the logics $l = (\Sigma, \vdash)$ with $\Sigma$ finitely presentable in $\Ss$ and $\vdash$ is a compact consequence relation in $Cons_\Sigma$.
\qdr\efa

\bre\ \label{As-re} In the sequence of works, \cite{AFLM1}, \cite{AFLM2}, \cite{AFLM3} is proven that  the category $\cA_s$ of Blok-Pigozzi algebraizable logics (\cite{BP}) and $\Ls$-morphisms that induces algebraizing pairs preserving functions on the formula algebras is a relatively complete $\omega$-accessible category (\cite{AR}).
\qdr\ere

\bre\ \label{Lsdefect} The \und{fundamental defect} of $\Ls$ is that the presentations of classical logic, for instante in the signatures $\Sigma = (\neg, \ra)$ and $\Sigma' = (\neg', \vee')$, {\em are not $\Ls$-isomorphic}: this deficiency is inherited from $\Ss$, because the $\Ss$-morphims are too strict.
\qdr\ere

\subsection{New results on categories of signatures and of logics}

\hp  The category $\Lf$ is the category of propositional logics and {\em flexible} translations
as morphisms.  This is a category ``built above'' the category $\Sf$, that
is, there is an obvious forgetful functor $U_f : \Lf \ \lra \ \Sf$. The categories $\Sf$ and $\Lf$  are considered in the literature (\cite{BC}, \cite{BCC1}, \cite{BCC2}, \cite{CG}, \cite{FC}), but  with a different the emphasis: Here we provide a more systematic analysis of category $\Lf$ and its relation with $\Sf$ and $\Ls$.

The objects of $\Lf$ are logics $l= (\Sigma, \vdash)$, as described in subsection 2.1.

If $l=(\Sigma,\vdash), l'=(\Sigma',\vdash')$ are logics then a
\emph{flexible translation morphism} $h : l \lra l'$ in $\Lf$ is a {\em flexible} signature morphism
$f : \Sigma \lra \Sigma'$ in $\Sf$ such that ``preserves the consequence relation'', that
is, for all $\Gamma \cup\{\psi\} \sub F(\Sigma)$, if $\Gamma \vdash \psi$
then $\overset{\vee}\to{h}[\Gamma] \vdash' \overset{\vee}\to{h}(\psi)$.  Composition and
identities are similar to $\Sf$.

As in $\Ls$,  $\Lf$ has  natural notions of direct and inverse image logics under a $\Sf$-morphism (just replace $\hat{f}$ by $\check{f}$) and they have good properties. 
For instance:

\bfa \label{Lfimage-fa}%
 Let $f : \Sigma \lra
\Sigma'$ be a $\Sf$-morphism and let $l=(\Sigma,\vdash),
l'=(\Sigma',\vdash')$ be logics $l,l' \in Obj(\Lf)$. Then
$(i)^{\star}$ and $(i)_{\star}$ hold, and $(ii)^{\star}, (ii)_{m}$
and $(ii)_{\star}$ are equivalent:
\item[$(i)^{\star}$] if $l'=(\Sigma',\vdash') \in Obj(\Lf)$ then
    $f^{\star}(l') = (\Sigma, \vdash_{f^{\star}(\vdash')}) \in Obj(\Lf)$;
   \item[$(i)_{\star}$] if $l=(\Sigma,\vdash) \in Obj(\Lf)$ then
    $f_{\star}(l) = (\Sigma', \vdash_{f_{\star}(\vdash)}) \in
    Obj(\Lf)$.
    \item[$(ii)^{\star}$] $\vdash\ \leq\ f^{\star}(\vdash')$;
    \item[$(ii)_{m}$] $f : (\Sigma,\vdash) \lra (\Sigma',\vdash')$ is $\Lf$-morphism;
    \item[$(ii)_{\star}$] $f_{\star}(\vdash)\ \leq'\ \vdash'$.
\qdr\efa

\bre\ \label{Lfadjunts-re}%
It follows easily from the facts above that the forgetful functor
$U_f : \Lf \lra \Sf$  : $((\Sigma, \vdash)$ $\overset{h}\to\lra (\Sigma',\vdash'))
\mapsto (\Sigma \overset{h}\to\lra \Sigma')$ has left and right adjoint
functors: the left adjoint  $\bot_f : \Ss \lra \Ls$ and the right adjoint $\top_f :
\Ss \lra \Ls$ take a signature $\Sigma$ to, respectively, $\bot_f(\Sigma) =
(\Sigma, \vdash_{min})$ (the first element of $Cons_{\Sigma}$) and
$\top_f(\Sigma) = (\Sigma, \vdash_{max})$ (the last element of $Cons_{\Sigma}$). Moreover, $U_f \circ \bot_f = Id_{\Sf} = U_f \circ \top_f$ and  $U_f$ preserves all limits and colimits that exists in $\Sf$.
\qdr\ere


\bre\ \label{Lflim} It is known that $\Lf$ has weak products, coproducts and some pushouts, and in the Remark above we see that $U_f$ preserves limits and colimits. As $U_f$ also "lift" limits and colimits -- the constructions in $\Lf$ are analogous to in $\Ls$
, just replace $\hat{f}$ by $\check{f}$ -- then given a small category $\cI$, $\Lf$ is $\cI$-complete (respectively, $\cI$-cocomplete) {\em if and only if} $\Sf$ is $\cI$-complete (respectively, $\cI$-cocomplete). Thus the Corollary \ref{Sfcolim-co} entails that $\Lf$ has colimits for any (small) diagram "in $\Ls$" (i.e., obtained via $(+) : \Ss \ \lra \ \Sf$), in particular, it has all unconstrained fibrings (= coproducts) and the constrained fibrings (= pushouts)  "based in $\Ls$".
\qdr\ere



\bre\ The fact of the formula algebra functions induced by  $\Sf$-morphisms "increase complexity" (see  Proposition \ref{reg-pr}.(a) for the precise statement) impose many limitations on $\Lf$. For instance:\\
(a) In \cite{CG} is shown that $\Lf$ solves the identity problem for the presentations of classical logic in terms of the (weaker) concept of {\em equipollence of logics}\footnote{We thank professor Marcelo Coniglio for that reference.}. But $\Lf$ does not solve problem of identity for the  presentations of classical logic  in terms of $\Lf$-isomorphisms. \\
(b) $\Lf$ has  weak terminal object but does not have terminal object; analogous statements holds  in general for $\Lf$ concerning (weak) products. 
\qdr\ere

The results below, together with  
 \ref{Lfadjunts-re}, \ref{Lflim}, constitute a strong evidence that the all defects in $\Lf$ are inherited from $\Sf$.

\bte\ \label{adjun-te} The signature adjunction $\Ss \overset{(+)_S}\to{\underset{(-)_S}\to\rlas} \Sf$ "lifts", via the forgetful functors $U_s$ and $U_f$, to  a logic adjunction $\Ls \overset{(+)_L}\to{\underset{(-)_L}\to\rlas} \Lf$, i.e.: \\ 
$\bullet$ \ $U_f \circ (+)_L  = (+)_S \circ U_s$; \\
$\bullet$ \ $U_s \circ (-)_L  = (-)_S \circ U_f$;\\
$\bullet$ \ $U_s \eta_L  =  \eta_S U_s$;\\
$\bullet$ \ $U_f \varepsilon_L   = \varepsilon_S U_f $.
\\
 Moreover, the following relations hold:\\
$\bullet$ \ $(+)_L \circ \bot_s = \bot_f \circ(+)_S$;\\
$\bullet$ \ $(+)_L \circ \top_s = \top_f \circ(+)_S$;\\
$\bullet$ \ $(-)_L \circ \bot_f = \bot_s \circ(-)_S$;\\
$\bullet$ \ $(-)_L \circ \top_f  \leqs \top_s \circ(-)_S$.
\ete
\bdm We provide only (in this moment), the definitions of the (faithful) functors: \\ 
\hfl $(+)_L : \Ls \ \lra \ \Lf$ \ : \hfl \\ 
\hfl $((\Sigma, \vdash) \overset{f}\to\lra\ (\Sigma', \vdash'))$ \ $\mapsto$ \ $((\Sigma, \vdash ) \overset{(j_{\Sigma'} circ f)^\flat}\to\lra\ (\Sigma',\vdash'))$; \hfl\\
\hfl $(-)_L : \Lf \ \lra \ \Ls$ \ : \hfl \\ \hfl $((\Sigma, \vdash) \overset{h}\to\lra\ (\Sigma', \vdash'))$ \ $\mapsto$ \ $(((F(\Sigma)[n])_{n \in \N},(j_{\Sigma})_\star(\vdash)) \overset{(\check{h}\rest_n)_{n\in\N}}\to\lra\ ((F(\Sigma' )[n])_{n \in \N},(j_{\Sigma'})_\star (\vdash')))$.\hfl
\qdr\edm


{
\begin{picture}(120,100)
\setlength{\unitlength}{.5\unitlength}
\thicklines
\put(32,20){$\Ss$}
\put(32,162){$\Ls$}
\put(09,39){\vector(0,1){115}}
\put(39,154){\vector(0,-1){115}}
\put(69,39){\vector(0,1){115}}
\put(75,177){\vector(1,0){380}}
\put(455,157){\vector(-1,0){380}}
\put(487,162){$\Lf$}
\put(15,90){{\small $\bot_s$}}
\put(45,90){{\small $U_s$}}
\put(75,90){{\small $\top_s$}}
\put(245,185){$(+)_L$}
\put(245,135){$(-)_L$}
\put(487,20){$\Sf$}
\put(435,90){{\small $\top_f$}}
\put(465,90){{\small $U_f$}}
\put(495,90){{\small $\bot_f$}}
\put(245,-5){$(+)_S$}
\put(245,45){$(-)_S$}
\put(462,39){\vector(0,1){115}}
\put(492,154){\vector(0,-1){115}}
\put(522,39){\vector(0,1){115}}
\put(455,35){\vector(-1,0){380}}
\put(75,15){\vector(1,0){380}}
\end{picture}}

\bte \label{LfKleisli-te} The signature monad  $\cT_S = (T_S, \eta_S, \mu_S)$  associated to the signature adjunction $ (\eta_S, \varepsilon_S) : \Ss \overset{(+)_S}\to{\underset{(-)_S}\to\rlas} \Sf$ (i.e., $\mu_S = (-)_S\varepsilon_S(+)_S$) "lifts" to a logic monad $\cT_L = (T_L, \eta_L, \mu_L)$  associated to the signature adjunction $ (\eta_L, \varepsilon_L) : \Ls \overset{(+)_L}\to{\underset{(-)_L}\to\rlas} \Lf$ (i.e., $\mu_L = (-)_L\varepsilon_L(+)_L$) and is such that $Kleisli(\cT_L) = \Lf$. Moreover, the functors $(+)_L$ and $(-)_L$ are precisely the canonical functors associated to the adjunction  of the Kleisli category of a monad. 
 \qdr\ete

\subsection{The appropriate categories of logics}


\hp In \ref{Lsdefect} we saw that the fundamental defect of logical category  $\Ls$ is due to the strictness of $\Ss$-morphisms and, analogously, in the previous subsection we saw that the deficiencies of the logical category $\Lf$ are inherited from the signature category $\Sf$.  Here we introduce  new categories of logics, still with logics as objects, one of them  satisfies {\em simultaneously} all the four natural requirements described in the Introduction.

\bct\ We will write $\Qf$ for the {\bf q}uotient category of $\Lf$ by the relation of interdemonstrability: the objects of $\Qf$ are the logics $l = (\Sigma, \vdash)$ and $\Qf((\Sigma, \vdash), (\Sigma', \vdash')) \ := \ \{ [f] : f \in \Lf((\Sigma, \vdash), (\Sigma', \vdash'))\}$, where $[f] := \{ g \in \Lf((\Sigma, \vdash), (\Sigma', \vdash')) : f \sim g \}$ and $f \sim g$ iff $(\overset{\vee}\to{f}/\!\!\dashv \vdash)  \ = \ (\overset{\vee}\to{g}/\!\!\dashv  \vdash)  \ : \  F(\Sigma)/\!\!\dashv  \vdash \ \lra \  F(\Sigma')/\!\! \dashv'  \vdash' \}$. Clearly, the relation $\sim$ is a congruence relation in the category $\Lf$ and we can take $\Qf := \Lf/\!\!\sim$.
\qdr\ect

\bfa\ (a) \ $\Qf$ has terminal object, coequalizers, {\em weak} products and {\em weak}  coproducts. \\
(b) \ The problem of identity for the  presentations of classical logic (without logical constants) is solved in terms of $\Qf$-isomorphisms.\\
(c) \ Any presentation of classical logic $l$ are $\Qf$-strongly rigid, i.e., if $f : l \ra l$ is a $\Lf$-morphim, then $[f] = [id] \in \Qf(l,l)$.\\
(d) By Proposition 4.3 in \cite{CG}, the logics $l$, $l'$ are \und{equipollent} iff they are $\Qf$-isomorphic.
\qdr\efa

\bct\ A logic $(\Sigma, \vdash)$ is {\em congruential} if, for each $c_{n} \in \Sigma_{n}$ and each $\{(\varphi_0,\psi_0), \ldots, (\varphi_{n-1}, \psi_{n-1}) \}$ such that $\varphi_0 \dashv \vdash \psi_0, \ldots, \varphi_{n-1} \dashv \vdash  \psi_{n-1}$, then $c_{n}(\varphi_0, \ldots, \varphi_{n-1}) \dashv \vdash c_{n}(\psi_0, \ldots, \psi_{n-1})$. It follows that if $\vartheta_{0}, \vartheta_{1} \in F(\Sigma)$ are such that $ var(\vartheta_0) = var(\vartheta_1) = \{x_{i_{0}}, \ldots, x_{i_{n-1}}\}$ and $\vartheta_0\dashv \vdash \vartheta_1$ then $\vartheta_0[\vec{x}\mid \vec{\varphi}] \dashv \vdash \vartheta_1[\vec{x}\mid \vec{\psi}]$. Clearly, the presentations of classical logic are congruential logics.
\qdr\ect
 
\bct\ Denote $\Lc$ the full subcategory of $\Lf$ whose objects are the congruential logics. 
\qdr\ect

\bpr\ $\Lc$ is a reflective subcategory of $\Lf$, i.e. $i : \Lc \hookr \Lf$ has a left adjoint $c: \Lf \ra \Lc$, moreover, the underlying signatures of the logics $l$ and $c(l)$ coincide. 
\epr

\bct \ Analogously to \ref{As-re}, denote $\cA_f$, the category of Blok-Pigozzi algebraizable logics and $\Lf$-morphisms that preserves algebraizable pairs (well defined). Let  $Lind(\cA_f) \sub \cA_f$, the full subcategory of  Lindenbaum algebrizable logics, i.e.  an algebraizable logic $l = (\Sigma, \vdash)$ is {\em Lindenbaum algebraizable} iff for each $\varphi, \psi \in F(\Sigma)$, \
$ \varphi \dashv  \vdash \psi$ $\Lra$
        $\vdash \varphi\Delta\psi$ (well defined).\\
\qdr\ect
        
\bfa \ (a) \ $Lind(\cA_f) \sub \Lc$.\\
(b) \ $Lind(\cA_f) \hookr \cA_f$  is a reflective  subcategory.
\qdr\efa

\bct\ It follows from the Proposition above that $\Lc$ has coproducts: it is the  "congruential closure" of the coproduct in $\Lf$ of a discrete diagram in $\Lc$.
\qdr\ect

\bct\ Let $\Qc$ denote the full subcategory of $\Qf$ whose objects are the congruential logics. Remark that $\Qc$ coincide with $Q(\Lc)$, the quotient of the (sub)category $\Lc$. 
\qdr\ect

{\em We propose that $\Qc$ is a convenient category to perform combinations of logics. The rest of the subsection is devoted to justify this claim.}

\bfa If $i : \Qc \hookr \Qf$ is the inclusion functor, then $i$  has a left adjoint $c: \Qf \ra \Qc$. 
\qdr\efa


{
\begin{picture}(120,100)
\setlength{\unitlength}{.5\unitlength}
\thicklines
\put(32,20){$\Qf$}
\put(32,162){$\Lf$}
\put(39,154){\vector(0,-1){115}}
\put(75,177){\vector(1,0){380}}
\put(455,157){\vector(-1,0){380}}
\put(487,162){$\Lc$}
\put(45,90){{\small $q_f$}}
\put(245,185){$c$}
\put(245,135){$i$}
\put(487,20){$\Qc$}
\put(465,90){{\small $q_f^c$}}
\put(245,-5){$c$}
\put(245,45){$i$}
\put(492,154){\vector(0,-1){115}}
\put(455,35){\vector(-1,0){380}}
\put(75,15){\vector(1,0){380}}
\end{picture}}

\bre\
(a) \ The proposition above, gives us a natural interpretation of the "ubiquity" of the logical systems considered. \\
(b) \ As in $\Qf$, the problem of identity for the  presentations of classical logic (without logical constants) is solved in terms of $\Qc$-isomorphisms.\\
(c) \ {\em Problem:} To describe explicitly congruential closures of not well-behaved logics as paraconsistent logics.
\qdr\ere

\bpr\ 
Given $f \in \Lc(l,l')$, then $[f]$ is a $\Qc$-isomorphism iff $f$ is an "weak  equivalence", i.e., it holds:\\
$\bullet$ \ {\em conservative translation:} \ $\all \Gamma \cup \{\psi\} \sub F(\Sigma)$ \ $\Gamma \vdash \psi$ \ $\Lra$ \ $\check{f}[\Gamma] \vdash' \check{f}(\psi)$ ;\\
$\bullet$ \ {\em denseness}: \ $\all n \in \N$\ $\all \vartheta' \in F(\Sigma')[n]$\ $\exists \vartheta \in F(\Sigma)[n]$ \ $\vartheta' \dashv ' \vdash \check{f}(\vartheta)$. 
\qdr\epr

\bct\ The {\em coproducts} in $\Qc$ are obtained taking first a  cone coproduct in $\Lf$:  the vertex in $\Qc$ is the congruential closure of the vertex in $\Lf$ and the cocone arrows in $\Qc$ are the classes of equivalence of the cocone arrows in $\Lf$  (the congruential property is decisive in proof of uniqueness). The {\em coequalizers} in $\Qc$ are obtained taking first a  cone coequalizer in $\Qf$ and then taking the induced cone in $\Qc$ obtained by the reflection functor $c: \Qf \ra \Qc$.
\qdr\ect

\bpr\ $\Qc$ is a cocomplete category. 
\qdr\epr

\bct\ A  congruential logics $(\Sigma, \vdash)$ is of {\em finite type} if it has a finite support signature ($card(\bigcup_{n \in \N} \Sigma_n) < \omega$) and is the congruential closure of a consequence relation over $\Sigma$ that is  generated by substitutions with a finite set of axioms and a finite set of (finitary) inference rules.
\qdr\ect

\bfa\ There is only an $2^{2^{\aleph_0}}$ set of classes of $\Qc$-isomorphism of finite type congruential logics. 
\qdr\efa 
 
\bpr\ Any congruential logic is a colimit in $\Qc$ of a directed diagram of congruential logics of finite type. 
\qdr\epr

\bpr\ 
 Let $(I, \leqs)$ be an upward directed poset and let $D : (I, \leqs) \ra \Qc$, $i \ra j$ \ $\mapsto$ \ $(\Sigma_i, \vdash_i) \overset{[h_{ij}]}\to\ra (\Sigma_j, \vdash_j)$ be a diagram. Then $colim (D) = ((\Sigma_j, \vdash_j) \overset{[\alpha^{+}_{j}]}\to\ra (\sqcup_{i \in I} \Sigma_i, \vdash)_{j\in I}$,where $\alpha_j \in \Ss(\Sigma_j, \sqcup_{i \in I} \Sigma_i)$ and for each $\Gamma \cup \{\varphi\}\sub F(\sqcup_{i \in I} \Sigma_i)$, \
$\Gamma \vdash \varphi$ \ $\Lra$ \\ $\exists \Gamma' \sub_{fin} \Gamma$, \ $\exists j \geqs I'= \{i \in I:$ for some $n \in \N$  some $n$-ary connective $(c_n,i) \in \sqcup_{i \in I} \Sigma_i$ occurs in $ \psi \in  \Gamma' \cup \{\varphi\}\} \sub_{fin} I$, \
$\Gamma'^{(j)} \vdash_j \varphi^{(j)}$, where
 for $\psi \in  \Gamma' \cup \{\varphi\}$ and each subformula $\theta$ of $\psi$: \\
$\bullet$ \  $\theta^{(j)} = x_l$,  if $\theta = x_l$;\\
$\bullet$ \ $\theta^{(j)} = \check{h}_{ij}((c_n,i))[x_0|\gamma_0^{(j)}, \ldots, x_{n-1}|\gamma_{n-1}^{(j)}]$, if $\theta = (c_n,i)(\gamma_0, \ldots, \gamma_{n-1})$.
\qdr\epr

\bpr\  In $\Qc$, a congruential logic of finite type is finitely presentable. Thus 
a  congruential logic is {\em finitely presentable} iff  it is a retract of  a congruential logic  of finite type.
\qdr\epr

\bte\ $\Qc$ is a finitely locally presentable category, i.e. it is a finitely accessible category complete/cocomplete.
\qdr\ete

A fundamental result in the theory of accessible categories (see \cite{AR}) ensures that an accessible category is complete iff its its cocomplete. But we will not  provide here explicit description of all limits in $\Qc$!


\bct\
A natural notion  of (Lindenbaum) algebrized logic is given by the triples $(\Sigma, \vdash, \Delta /\!\!\dashv \vdash)$ where $l = (\Sigma, \vdash)$ is a logic and
$\Delta \ \sub_{fin} \ F(\Sigma)[2]$ is a set of "equivalence formulas in the Lindenbaum sense" i.e.:\\
 { (a) } \ $\vdash \varphi \Delta \varphi$ \footnote{That is, if $\Delta = \{ \Delta_u : u <v\}$, then
        $\vdash \varphi \Delta_{u} \varphi$, for all $u<v$.}; \\
{ (b)} \ $\varphi\Delta\psi \vdash \psi\Delta\varphi$;\\
{ (c)} \ $\varphi\Delta\psi, \psi\Delta\vartheta \vdash
        \varphi\Delta\vartheta$;\\
{ (d)} \ $\varphi_{0}\Delta\psi_{0}, \ldots,
        \varphi_{n-1}\Delta\psi_{n-1} \vdash
        c_n(\varphi_{0},\ldots,\varphi_{n-1})\Delta
        c_n(\psi_{0},\ldots,\psi_{n-1})$;\\ 
{ (e)} \ $ \varphi \dashv  \vdash \psi$ iff
        $\vdash \varphi\Delta\psi$.\\
  Clearly, the underlying logic of $(l, \Delta /\!\!\dashv  \vdash)$ is congruential.
\qdr \ect
 
\bct\ The corresponding category $\cA$ of algebrized logics has as morphisms $f : (l, \Delta /\!\!\dashv  \vdash) \ \lra \ (l',\Delta' /\!\!\dashv'  \vdash')$ the 
$\Delta'  \vdash' \overset{\vee}\to{f} [\Delta]$ (or $\Delta' \dashv' \vdash \overset{\vee}\to{f} [\Delta]$). Composition and identities are as in $\Qc$
qdr\ect

\vspace{0.5cm}
We finish the section with the following diagram:

\begin{picture}(170,170)
\setlength{\unitlength}{.9\unitlength} \thicklines 
\put(10,20){$Lind(\cA_f)$}
\put(32,162){$\cA_f$} \put(35,154){\vector(0,-1){115}}
\put(45,39){\vector(0,1){115}}
\put(180,154){\vector(0,-1){115}}
\put(190,39){\vector(0,1){115}}
\put(325,154){\vector(0,-1){115}}
\put(335,39){\vector(0,1){115}}
 \put(65,168){\vector(1,0){100}}
  \put(65,25){\vector(1,0){100}}
  \put(177,162){$\Lf$} 
  \put(20,90){$L$} \put(50,90){$j$}
\put(100,175){\small $incl$} 
\put(100,10){\small $incl$} 
\put(177,20){$\Lc$} \put(320,20){$\Qc$} \put(320,162){$\Qf$}\put(195,90){$i$}
\put(165,90){$c$} \put(340,90){$\bar{i}$} \put(310,90){$\bar{c}$} \put(250,8){$q^c$} \put(250,175){$q$}
\put(210,168){\vector(1,0){100}} \put(210,25){\vector(1,0){100}}
\end{picture}

\section{Final remarks and future works}

\bct\ It should be remarked that all the categories of logics considered have the same  objects (of combinatorial nature) $l = (\Sigma, \vdash)$ and the morphisms in all that categories area, in some sense,  translation morphisms: thus we believe that the new categories introduced ($\Lc, Lind(\cA_f), \Qf, \Qc$) constitute a natural and simple solution of the deficiencies on $\Ls$ and $\Lf$, that satisfies the natural requirements (i),...,(iv) in the Introduction.
\qdr \ect


\bct\ In the definition of the quotient category $\Qf$ (respec. $\Qc$), is considered a congruence on  category $\Lf$ (respec. $\Lc$), that is induced by a {\em pre-order} relation on the category: thus $\Qf$ (respec. $\Qc$) has a natural structure of category enriched by the category of posets (and increasing functions) and this should be explored in the future. 
\qdr \ect

\bct \ Another possible approach to overcome  the deficiencies of the categories of logics that are inherited from the defects of the underlying categories of 
\cite{AFLM2}, \cite{AFLM3}): A  possible notion candidate for this
task is concept of {\em operad}, as mentioned in the section 5 of \cite{AFLM3}. An  operad can be
understood as an axiomatization of the behavior of a collection of
finitary operations on a set and such that is closed by the
formation of the derived operations by composition. A candidate
for the category of signatures would be the category of all
operads: a category that somehow contains the categories proposed
in the two series of papers mentioneds, but that it has good
categorial properties (is complete, is cocomplete,...). Moreover
the category of logical systems built on this new category of
signatures would allow, in principle, the interdefinability of
connectives as now the operations may satisfy relations, something
that would promote a satisfactory treatment of the problem of
identity. 
\qdr\ect


\bct \ In \cite{MP} is in developement an alternative to overcome the problems by a mathematical device generically called "Morita equivalence of logics", a notion borrowed from Ring and Module Theories: 
In the representation theory of rings, the category of rings  is functorially encoded into the category of categories: a ring $R$ is encoded by the  category of (left/right) linear representation of $R$  (respec. $R-Mod$ , $Mod-R$). In the same vein, it is proposed a encoding of  a general propositional logic by a diagram of categories and functors given by the quasivarieties canonically associated to the algebraizable logics (in the sense of \cite{BP})  connected with the given propositional logic. 
\qdr \ect

\end{document}